\documentclass[reqno,12pt]{amsart}
\usepackage{bbm}
\usepackage{amsmath}
\usepackage{txfonts}
\usepackage{stmaryrd}
\usepackage{amssymb}
\usepackage{amsfonts}
\usepackage{times}
\usepackage{mathrsfs}

\usepackage{amsthm}
\usepackage{enumerate}

\usepackage{framed}
\usepackage{lipsum}
\usepackage{color}

\def\CC{ \mathbb{C}}
\newcommand{\DD}{{\mathbb D}}

\def\R{\mathcal{R}}
\def\I{\mathcal{I}}
\def\D{\mathbb{D}}
\def\N{\mathbb N}

\def\hD{\hat{\mathcal{D}}}

\def\vp{\varphi}
\def\om{\omega}

\def\p{{\prime}}

\def\msk{\medskip}
\def\ol{\overline}

\begin{document}
\title[  ]{ Weighted  composition operators on weighted  Bergman   spaces induced by double weights}
\author{ Juntao Du,  Songxiao Li$\dagger$ and Yecheng Shi  }
\address{Juntao Du\\ Faculty of Information Technology, Macau University of Science and Technology, Avenida Wai Long, Taipa, Macau.}
\email{jtdu007@163.com  }
\address{Songxiao Li\\ Institute of Fundamental and Frontier Sciences, University of Electronic Science and Technology of China,
610054, Chengdu, Sichuan, P.R. China\newline
Institute of Systems Engineering, Macau University of Science and Technology, Avenida Wai Long, Taipa, Macau. } \email{jyulsx@163.com}

\address{Yecheng Shi\\ School of Mathematics and Statistics, Lingnan Normal University, Zhanjiang 524048, Guangdong, P. R. China }\email{ 09ycshi@sina.cn}

\subjclass[2000]{30H10, 47B33 }
\begin{abstract}  In this paper, we investigate the boundedness, compactness, essential norm  and the Schatten class of weighted composition operators $uC_\vp$ on  Bergman type spaces $A_\om^p $ with double weight $\om$. Let $X=\{u\in H(\D): uC_\vp:A_\om^p\to A_\om^p \mbox{ is bounded}\}.$
For some regular weights $\om$, we obtain that $X=H^\infty$   if and only if $\vp$ is a finite Blaschke product.
\thanks{$\dagger$ Corresponding author.}
\vskip 3mm \noindent{\it Keywords}: Weighted composition operator, weighted Bergman space,  double   weight.
\end{abstract}
 \maketitle

\section{Introduction}
Let $\D$ be the open unit disk in the complex plane, and $H(\D)$ the class of all functions analytic on $\D$. Let $\vp$ be an analytic self-map of $\DD$ and $u\in H(\D)$. The weighted  composition operator, denoted by $uC_\vp$, is defined on $H(\D)$ by
$$(uC_\vp f)(z)=u(z)f(\vp(z)), ~~~~~f\in H(\D).$$

 For $0<p<\infty$, $H^p$ denotes the  Hardy space, which consisting  of all functions  $f\in H(\D)$ satisfied
$$
\|f\|_{H^p}^p=\sup_{0<r<1} \frac{1}{2\pi}\int_0^{2\pi}|f(re^{i\theta})|^p d\theta<\infty.
$$
As usual, $H^\infty$ is the set of bounded analytic functions in $\D$.

We say that $\mu$ is a weight, when $\mu$ is  radial and  positive  on $\D$. Suppose that $\om$ is an integrable weight  on $(0,1)$.  Let $\hat{\om}(r)=\int_r^1 \om(s)ds$ for $r\in(0,1)$.  We say that $\om$ is regular, denoted by $\om\in\R$, if there is a constant $C>0$ depending on $\om$, such that
$$\frac{1}{C}<\frac{\hat{\om}(r)}{(1-r)\om(r)}<C,\,\,\mbox{ when } 0<r<1.$$
We say that $\om$ is rapidly increasing, denoted by $\om\in\I$, if
$$\lim_{r\to 1} \frac{\hat{\om}(r)}{(1-r)\om(r)}=\infty.$$
Let
$$v_{\alpha,\beta}(r)=(1-r)^\alpha\left(\log \frac{e}{1-r}\right)^\beta.$$
After a calculation, we have the following  typical examples of regular and rapidly increasing weights, see  \cite{PjaRj2014book}, for example.
\begin{enumerate}[(i)]
  \item When $\alpha>-1$ and $\beta\in \mathbb{R}$, $v_{\alpha,\beta}\in\R$;
  \item When $\alpha=-1$ and $\beta<-1$, $v_{\alpha,\beta}\in\I$  and $\left|\sin\left(\log\frac{1}{1-r}\right)\right|v_{\alpha,\beta}(r)+1\in\I$.
\end{enumerate}

In \cite{Pja2015}, Pel\'aez introduced the set of  double weights, denoted by $\hD$, which includes $\I\cup\R$.
 We say that $\om\in \hD$ if there is a constant $C>0$ such that
$\hat{\om}(r)<C\hat{\om}(\frac{1+r}{2}) ,\,\,\mbox{ when } 0<r<1.$ We should remark that  the most part of the results in \cite{PjaRj2014book}, which presented in the context of regular and rapidly increasing weights, continue to hold for the wider class $\hD$.
More details about $\I,\R$ and $\hD$ can be seen in \cite{Pja2015,PjaRj2014book,PjaRj2016}.

For $0<p<\infty$ and $\om\in\hD$, the weighted Bergman space $A_\om^p$ is the space of $f\in H(\D)$ for which
$$\|f\|_{A_\om^p}^p=\int_{\D}|f(z)|^p\om(z)dA(z)<\infty,$$
where $dA(z)=\frac{1}{\pi}dxdy$ is the normalized Lebesgue area measure on $\D$. When $\om(t)=(1-t)^\alpha(\alpha>-1)$, the space $A_\om^p$ becomes the classical weighted Bergman space $A_\alpha^p$. For classical Bergman space $A_\alpha^p$, we refer \cite{CccMbd1995,DpSa2004book,Zk2007book} and references therein.
In many respects, the Hardy space $H^p$ is the limit of $A_\alpha^p$ as $\alpha\to -1$. But it is a rough estimate since most of the finer function-theoretic properties of the classical weighted Bergman space $A_\alpha^p$ can not carry over to the Hardy space $H^p$. As we know,   $A_\om^p$ induced by regular weights have similar properties with  $A^p_\alpha$. But many results in \cite{Pja2015,PjaRj2014book,PjaRj2015,PjaRj2016,PjaRj2016jmpa,PjaRjSk2018jga}
show that spaces $A_\om^p$ induced by rapidly increasing weights, lie ``closer" to $H^p$ than any $A^p_\alpha$.

In  \cite{CzZr2004jlmc},   $\check{\mathrm{C}}$u$\check{\mathrm{c}}$kovi$\check{\mathrm{c}}$ and   Zhao
characterized the boundedness and compactness of weighted composition operators mapping on Bergman space $A^p_\alpha$ by using Berezin transform.
Three years later, they investigated  weighted composition operators between different Bergman spaces and Hardy spaces   in \cite{CzZr2007ijm}.
In \cite{PjaRj2016},   Pel\'aez and   R\"atty\"a characterized the Schatten class of Toeplitz operators induced by a positive Borel measure on $\D$ and the reproducing kernel of the Bergman space $A_\om^2$ when $\om\in\hD$.
In \cite{ZlHs2015ams},  Zhao and  Hou proved that,  for $A_\alpha^p$, the finite Blaschke product is the only composition symbol that the induced weighted composition operator is bounded if and only if the weighted symbol defines a bounded multiplication operator. The similar result for Hardy space $H^p$ can be seen in \cite{CmdHag2003}.

Motivated by \cite{CzZr2004jlmc,CzZr2007ijm,PjaRj2016}, under the assumption that $\om\in \hD $ and $\mu$ is a positive Borel measure, we investigate the boundedness, compactness and  essential norm of  $uC_\vp:A_\om^p\to L_\mu^q$ and the Schatten class of $uC_\vp:A_\om^2\to A_\om^2$.
Motivated by \cite{ZlHs2015ams}, we get that, for some $\om\in\R$,   $X=H^\infty$ if and only if $\vp$ is a finite Blaschke product. Here
$$X=\{u:u\in H(\D) \mbox { and } uC_\vp:A_\om^p\to A_\om^p \mbox{ is bounded}\}.$$

Throughout this paper, the letter $C$ will denote  constants and may differ from one occurrence to the other.
The notation $A \lesssim B$ means that there is a positive constant C such that $A\leq CB$.
The notation $A \approx B$ means $A\lesssim B$ and $B\lesssim A$.\msk

\section{Auxiliary results}

In this section we formulate and prove  several auxiliary results which will be used in the
proofs of  main results in this paper.\msk

{\noindent\bf Lemma 1.} {\it Assume that $\om\in\hD$, $r\in (0,1]$ and $\om_*(r)=\int_r^1 s\om(s)\log\frac{s}{r}ds$. Then  the following statements hold.
\begin{enumerate}[(i)]
  \item $\om_*\in\R$ and $\om_*(r)\approx (1-r)\hat{\om}(r)$ as $r\to 1$;
  \item There are  $1<a<b<+\infty$ and $\delta\in [0,1)$, such that
\begin{align}
\frac{\om_*(r)}{(1-r)^a} \;\; \mbox{is decreasing on}\;\; [\delta,1)\;
\mbox{and}\;\; \lim_{r\to 1}\frac{\om_*(r)}{(1-r)^a}=0;   \label{0521-1}  \\
\frac{\om_*(r)}{(1-r)^b} \;\; \mbox{is increasing on}\;\; [\delta,1)\;\;
\mbox{and}\;\; \lim_{r\to 1}\frac{\om_*(r)}{(1-r)^b}=\infty; \label{0521-2}
\end{align}
\item $\om_*(r)$ is decreasing on $[\delta,1)$ and $\lim\limits_{r\to 1} \om_*(r)=0.$
\end{enumerate} }

\begin{proof} By \cite[Lemmas A and 9]{PjaRj2016} and (1.19) in \cite{PjaRj2014book}, ({\it i}) and ({\it ii}) hold.
 ({\it iii}) follows by ({\it ii}) and $\om_*(r)=\frac{\om_*(r)}{(1-r)^a}(1-r)^a$.
\end{proof}

  {\noindent\bf Remark 1.} We observer that $z=0$ is the logarithmic singular point of $\om_*$.
So, for any fixed $r_0\in (0,1)$, we have $\om_*(r)\approx (1-r)\hat{\om}(r)$ for $r_0\leq r<1.$
For simplicity, suppose $\om_*$ and $\hat{\om}$ are radial, that is, $\om_*(z)=\om_*(|z|)$ and $\hat{\om}(z)=\hat{\om}(|z|)$ for all $z\in\D$.  \msk

Suppose $\mathbb{T}$ is the boundary of $\D$ and $I\subset\mathbb{T}$ is an interval. The Carleson square $S(I)$ can be defined as
$$S(I)=\{re^{it}:e^{it}\in I, 1-|I|\leq r<1\},$$
where $|I|$ denotes the Lebesgue measure of $I$.
For convenience, for each $a\in\D\backslash\{0\}$, we define
$$I_a=\left\{e^{i\theta}:|\arg(ae^{-i\theta})|\leq \frac{1-|a|}{2}\right\}$$
and denote $S(a)=S(I_a)$. By (26) in \cite{Pja2015}, when $\om\in\hD$, we have
\begin{align}\label{1017-3}
\om(S(a))\approx \om_*(a), \mbox{ for all } a\in\D \mbox{ and } |a|\geq \frac{1}{2}.
\end{align}

The  following lemma is a straight result of \cite[Lemma 3.1]{Pja2015}(or \cite[Lemma 2.4]{PjaRj2014book}).\msk

{\noindent\bf Lemma 2.} {\it Suppose $\om\in\hD$ and $0<p<\infty$. There exists $\gamma_0>0$, if $\gamma>\gamma_0$,
we have
$$|F_{a,p,\gamma}(z)|\approx \frac{1}{\om(S(a))^\frac{1}{p}},\,\,  \|F_{a,p,\gamma}\|_{A_\om^p}\approx 1, \,\,\mbox{ when } a\in\D,\,\,z\in S(a),$$
and
$$\lim_{|a|\to 1}\sup_{|z|\leq r} |F_{a,p,\gamma} (z)|=0,\,\,\mbox{ when } r\in(0,1).$$
Here and henceforth,
$$F_{a,p,\gamma} (z)=\left(\frac{1-|a|^2}{1-\ol{a}z}\right)^{\frac{\gamma+1}{p}}\frac{1}{(\om(S(a)))^\frac{1}{p}}.$$
}

For simplicity, in the rest of this paper, we always assume that $\gamma$ is large enough so that Lemma 2 holds when we mention the function $F_{a,p,\gamma}$.

For a given Banach space $X$ of analytic functions on $\D$, a positive Borel measure $\mu$ on $\D$ is called a $q-$Carleson measure
for $X$, if   the identity operator $Id: X\to L^q(\mu)$ is bounded.
By \cite[Theorem 3.3]{Pja2015}, when $\om\in\hD$,  a Borel measure $\mu$ on $\D$ is a $q-$Carleson measure for $A_\om^p$  if and only if
$$\sup_{a\in\D}\frac{\mu(S(a))}{(\om(S(a)))^\frac{q}{p}}<\infty.$$
Moreover, $\|Id\|_{A_\om^p\to L^q_\mu}\approx \sup_{a\in\D}\frac{\mu(S(a))}{(\om(S(a)))^\frac{q}{p}}$.
Then we have the following lemma.\msk

{\noindent\bf Lemma 3.} {\it Suppose $\om\in \hD$, $\mu$ is a positive Borel measure on $\D$. Let $0<p\leq q<\infty$. For some (equivalently for all) large enough $\gamma$, we have
$$
\|Id\|_{A_\om^p\to L_\mu^q}^q \approx\sup_{a\in\D}\frac{\mu(S(a))}{\om(S(a))^\frac{q}{p}}
\approx\sup_{a\in\D}\int_{S(a)} |F_{a,p,\gamma}(z)|^q d\mu(z)\approx \sup_{a\in\D}\int_\D |F_{a,p,\gamma}(z)|^q d\mu(z).
$$
Here $Id$ is the identity operator.  }

\begin{proof} By Lemma 2, we have
\begin{equation}\label{0323-1}
\frac{\mu(S(a))}{\om(S(a))^{\frac{q}{p}}} \approx \int_{S(a)}|F_{a,p,\gamma}(z)|^q d\mu(z) \leq  \int_{\D}|F_{a,p,\gamma}(z)|^q d\mu(z), \mbox{ when } a\in\D.
\end{equation}
So,
$$\sup_{a\in\D}\frac{\mu(S(a))}{\om(S(a))^{\frac{q}{p}}}\approx \sup_{a\in\D}\int_{S(a)}|F_{a,p,\gamma}(z)|^q d\mu(z)
\lesssim  \sup_{a\in\D}\int_{\D}|F_{a,p,\gamma}(z)|^q d\mu(z).
$$
By Lemma 2 and  \cite[Theorem 3.3]{Pja2015} (also see  \cite[Theorem 2.1]{PjaRj2014book}), we obtain
$$\int_{\D}|F_{a,p,\gamma}(z)|^q d\mu(z)\lesssim \|Id\|_{A_\om^p\to L_\mu^q}^q \approx \sup_{a\in\D} \frac{\mu(S(a))}{\om(S(a))^{\frac{q}{p}}}.$$
 The proof is complete.
\end{proof}\msk

{\noindent \bf Lemma 4.} {\it  Suppose $0<p\leq q<\infty$, $\om\in\hD$, $\mu$ is a positive Borel measure on $\D$, and $\gamma$ is large enough.
Let $\frac{1}{2}<r<1$ and
$$N_r^*=\sup_{|a|>r}\int_\D |F_{a,p,\gamma}(z)|^q d\mu(z).$$
If $\mu$ is a $q$-Carleson measure for $A_\om^p$, then   $\mu_r=\mu|_{\D\backslash r\D}$ is also a $q$-Carleson measure for $A_\om^p$, where $r\D=\{z\in\D:|z|<r\}$.
Moreover, there is a $C>0 $, such that
\begin{equation}\label{0323-2}
\sup_{a\in\D} \frac{\mu_r(S(a))}{(\om(S(a)))^\frac{q}{P}}\leq CN_r^*.
\end{equation}}

\begin{proof}  It is obvious that $\mu_r$ is a $q$-Carleson measure for $A_\om^p$.
Let
$$N_r=\sup_{|a|\geq r}\frac{\mu(S(a))}{(\om(S(a)))^\frac{q}{p}}.$$
When $|a|\geq r$, $\frac{\mu_r(S(a))}{(\om(S(a)))^\frac{q}{p}}<N_r$ is obvious. When $|a|<r$, letting $k=\mbox{int}(\frac{1-|a|}{1-r})+1$, there exists $a_1,a_2,\cdots,a_k\in\D$ such that $ S(a)\cap {\D\backslash r\D} \subset \cup_{i=1}^k S(a_i) $ and $|a_i|=r$ for $i=1,2,\cdots,k$.
By Lemma 1, we have

\begin{align*}
\mu_r(S(a))
&\leq \sum_{i=1}^k \mu(S(a_i))      \leq N_r \sum_{i=1}^k (\om(S(a_i)))^\frac{q}{p}  \\
&\lesssim  N_r\left(\frac{1-|a|}{1-r}+1\right)\om_*(r)^\frac{q}{p}   \\
&\approx N_r\left(\frac{1-|a|}{1-r}+1\right)(1-r)^\frac{q}{p}\hat{\om}(r)^\frac{q}{p}   \\
&\leq N_r\left(\left(\frac{1-r}{1-|a|}\right)^{\frac{q}{p}-1}+\left(\frac{1-r}{1-|a|}\right)^{\frac{q}{p}}\right)
(1-|a|)^\frac{q}{p}\hat{\om}(a)^\frac{q}{p}  \\
&\lesssim  N_r  \om(S(a))^\frac{q}{p}.
\end{align*}
So, there exists $C >0$, such that
$$\frac{\mu_r(S(a))}{\om(S(a))^\frac{q}{p}}\leq C N_r .$$
By (\ref{0323-1}), we have $N_r \lesssim N_r^*.$ Therefore,  (\ref{0323-2}) holds.  The proof is complete.
\end{proof}

The following lemma can be proved in a standard way (see, for example,
Theorem 3.11 in \cite{CccMbd1995}). \msk

  \noindent{\bf Lemma 5.}  {\it  Suppose $0<p,q<\infty$, $\om\in\hD$ and $\mu$ is a  positive Borel measure. If\,     $T:A_\om^p\rightarrow L_\mu^q$ is linear and bounded,  then $T$ is  compact  if and only if whenever  $\{f_k\}$   is bounded in $A_\om^p$ and $f_k \rightarrow 0$  uniformly  on compact subsets of $\D$,   $\lim\limits_{k\to \infty}\|T f_{k }\|_{L_\mu^q}=0$.} \msk

The following lemma can be found in \cite{ZlHs2015ams} without a proof. For the benefits of the readers, we will prove it. \msk

\noindent{\bf Lemma 6. }{\it Let $\vp$ be an analytic self-map of $\D$. Then $\vp$ is a finite Blaschke product if and only if  $\lim\limits_{|w|\to 1} |\vp(w)|=1$.
}
\begin{proof} The sufficiency of the statement is obvious. Next we prove the necessity.

Suppose $\lim\limits_{|w|\to 1} |\vp(w)|=1$. Let $E\subset \D$ be compact. Then there exists a $r\in(0,1)$ such that $E\subset r\ol{\D}$ where $r\ol{\D}=\{z\in\D:|z|\leq r\}$.
Since $\lim\limits_{|w|\to 1} |\vp(w)|=1$, there is a $t\in (0,1)$ such that for all $|z|>t$, we have $|\vp(z)|>r$.
Therefore, $\vp^{-1}(E)\subset t\ol{\D}$.
By the continuity of $\vp$, $\vp^{-1}(E)$ is closed. So, $\vp^{-1}(E)$ is compact.
By the subsection 7.1.3 of  \cite{Rw2003}, $\vp$ is proper. By the subsection 7.3.1 of \cite{Rw2003}, $\vp$ is a finite Blaschke product.
The proof is complete.
\end{proof}\msk

 \section{ Main results and proofs}

 {\noindent\bf Theorem 1.}  {\it Assume  $\om\in\hD,0<p\leq q<\infty$, $u:\D\to\CC$ is a measurable function, $\vp$ is an analytic self-map of $\D$, and $\mu$ is a positive Borel measure on $\D$. Then
 $$\|uC_\vp\|_{A_\om^p\to L^q_\mu}^q  \approx  \sup_{a\in\D} \int_\D |F_{a,p,\gamma}(\vp(z))|^q  |u(z)|^q d\mu(z).$$ }

\begin{proof} By Lemma 2, we have
$$  \sup_{a\in\D} \int_\D |F_{a,p,\gamma}(\vp(z))|^q  |u(z)|^q d\mu(z)  \lesssim  \|uC_\vp\|_{A_\om^p\to L^q_\mu}^q .$$
 Let  $\upsilon(E)=\int_{\vp^{-1}(E)}|u(z)|^q d\mu(z)$ for all Borel set $E$. For all $f\in A_\om^p$, letting $w=\vp(z)$, by Lemma 3 we have
\begin{eqnarray}
\|uC_\vp f\|_{ L_\mu^q}^q &=&\int_\D |f(\vp(z))|^q |u(z)|^q d\mu(z)
=\int_\D |f(w)|^q d\upsilon(w)  \label{0313-1}\\
&=&\|f\|_{L_\upsilon^q}^q\leq\|Id\|_{A_\om^p\to L_\upsilon^q}^q \|f\|_{A_\om^p}^q\nonumber\\
&\lesssim & \|f\|_{A_\om^p}^q\sup_{a\in\D}\int_{\D}|F_{a,p,\gamma}(w)|^q d\upsilon(w)  \nonumber
\end{eqnarray}
Making the changing of variable $w=\vp(z)$, we have
$$ \int_{\D}|F_{a,p,\gamma}(w)|^q d\upsilon(w) =\int_\D |F_{a,p,\gamma}(\vp(z))|^q  |u(z)|^q d\mu(z). $$
The proof is complete.
\end{proof}

Let $X$ and $Y$ be  Banach spaces. Recall that the essential norm of linear operator $T:X\to Y$ is   defined by
$$\|T\|_{e,X\to Y}=\inf\{\|T-K\|_{X\to Y}: K \mbox{ is compact from }X \mbox{ to } Y\}.$$
Obviously $T:X\to Y$ is compact if and only if $\|T\|_{e,X\to Y}=0.$\msk

{\noindent\bf Theorem 2.}  {\it Assume  $\om\in\hD,1\leq p\leq q<\infty$, $u:\D\to\CC$ is a measurable function, $\vp$ is an analytic self-map of $\D$, and $\mu$ is a positive Borel measure on $\D$. If $uC_\vp:A_\om^p\to A_\mu^q$ is bounded, then
 $$\|uC_\vp\|_{e,A_\om^p\to L^q_\mu}^q  \approx  \limsup_{|a|\to 1} \int_\D |F_{a,p,\gamma}(\vp(z))|^q  |u(z)|^q d\mu(z).$$
}

\begin{proof} Since $uC_\vp:A_\om^p\to L_\mu^q$ is bounded, $u\in L_\mu^q$ and $\|u\|_{L_\mu^q}\leq \|uC_\vp\|_{A_\om^p\to L_\mu^q}$.

{\bf Upper estimate of $\|uC_\vp\|_{e,A_\om^p\to L_\mu^q}$. }

Suppose $f(z)=\sum_{k=0}^\infty \hat{f}_kz^k\in H(\DD)$. For $n\in\N$, let
$$
K_n f(z)=\left\{
\begin{array}{lc}
 \sum\limits_{k=0}^n \hat{f}_kz^k,& p>1; \\
  \sum\limits_{k=0}^n \left(1-\frac{k}{n+1}\right)\hat{f}_kz^k, & p=1,
\end{array}
\right.
$$
and
$R_n=Id-K_n.$

By \cite[Proposition 1 and Corollary 3]{Zk1991mmj}, when $1<p<\infty$, $K_n$ is bounded uniformly on $H^p$.
By \cite{Dr2011sm}, $\|K_n\|_{H^1\to H^1}\leq 1$.
So, when $p\geq 1$,  there is a $C=C(p)$ such that
\begin{eqnarray*}
\|K_n f\|_{A_\om^p}^p
&\leq& C\int_0^1 \om(s)sds\int_0^{2\pi}|f(se^{\mathrm{i}\theta})|^pd\theta \leq C\|f\|_{A_\om^p}^p,
\end{eqnarray*}
and
\begin{align}\label{1102-1}
\|R_n\|_{A_\om^p\to A_\om^p}=\|Id-K_n\|_{A_\om^p\to A_\om^p}\leq 1+\|K_n\|_{A_\om^p\to A_\om^p}\leq C+1 .
\end{align}
By Lemma 5 and Cauchy's estimate, $f\to \hat{f}_kz^k$ is compact on $A_\om^p$.
Therefore, $K_n:A_\om^p\to A_\om^p$ is compact. So, we have
\begin{align}\label{0314-1}
\|uC_\vp\|_{e,A_\om^p\to L_\mu^q}
= \|uC_\vp( K_n+R_n)\|_{e,A_\om^p\to L_\mu^q}
\leq \|uC_\vp R_n\|_{e,A_\om^p\to L_\mu^q}
\leq \|uC_\vp R_n\|_{A_\om^p\to L_\mu^q}.
\end{align}

For any fixed $r\in (0,1)$, by (\ref{0313-1}) we have
\begin{equation}\label{0314-2}
\|uC_\vp R_n f\|_{L_\mu^q}^q =\int_{\D\backslash r\D}|R_nf(w)|^qd\upsilon(w) + \int_{ r\D}|R_nf(w)|^qd\upsilon(w),
\end{equation}
where $\upsilon$ is defined in the proof of Theorem 1,
that is, for every Borel set $E\subset \D$, $\upsilon(E)=\int_{\vp^{-1}(E)}|u(z)|^q d\mu(z)$.

Let $\om_n=\int_0^1 s^n\om(s)ds$. Since $B_z^\om(\zeta)=\sum_{k=0}^\infty \frac{(\zeta\ol{z})^n}{2\om_{2n+1}}$ is the reproducing kernel of $A_\om^p$ (see \cite{PjaRj2016jmpa,PjaRjSk2018jga} for  example),   we have
\begin{eqnarray*}
|R_n f(w)|=\left|\langle R_n f, B_w^\om\rangle_{A_\om^2}\right|=\left|\langle  f, R_n B_w^\om\rangle_{A_\om^2}\right|
\lesssim \|f\|_{A_\om^p}  \|R_n B_w^\om\|_{H^\infty}.
\end{eqnarray*}
Here $\langle\cdot,\cdot\rangle_{A_\om^2}$ is the inner product induced by $\|\cdot\|_{A_\om^2}$
and $\|\cdot\|_{H^\infty}$ is the norm of the bounded analysis function space on $\D$.
When $|w|\leq r$, we have
$$
\|R_n B_w^\om\|_{H^\infty}\leq \frac{1}{n}\sum_{k=1}^\infty \frac{k r^{k-1}}{2\om_{2k+1}}  + \sum_{k=n+1}^\infty \frac{r^k}{2\om_{2k+1}}.
$$
By \cite[Lemma 6]{PjaRjSk2018jga}, $\sum\limits_{k=1}^\infty \frac{k r^{k-1}}{2\om_{2k+1}}$ is convergent and $\lim\limits_{n\to\infty}\sum\limits_{k=n+1}^\infty \frac{r^k}{2\om_{2k+1}}=0.$
So, for all $\varepsilon>0$, there is a $N=N(\varepsilon,\om,r)$, such that
$$|R_n B_w^\om|<\varepsilon , \mbox{ for all } |w|\leq r \mbox{ and } n>N.$$
So,  for all  $n>N$,
\begin{equation}\label{0314-3}
\int_{ r\D}|R_nf(w)|^qd\upsilon(w)\leq  \varepsilon^q \|u\|_{L_\mu^q}^q \|f\|_{A_\om^p}^q  .
\end{equation}
Let $\upsilon_r=\upsilon|_{\D\backslash r\D}$. By (\ref{1102-1}),  Lemmas 3 and 4,  we have
\begin{eqnarray}
\int_{\D\backslash r\D}|R_nf(w)|^qd\upsilon(w)
&=&\int_\D |R_n f(w)|^q d\upsilon_r(w)
   \lesssim \|R_n f\|_{A_\om^p}^q\sup_{a\in\D}\frac{\upsilon_r(S(a))}{\om(S(a))^\frac{q}{p}} \nonumber \\
&\lesssim&   \|R_n f\|_{A_\om^p}^q  \sup_{|a|>r}\int_\D |F_{a,p,\gamma}(\vp(z))|^q |u(z)|^q d\mu(z) \nonumber\\
&\lesssim& \| f\|_{A_\om^p}^q  \sup_{|a|>r}\int_\D |F_{a,p,\gamma}(\vp(z))|^q |u(z)|^q d\mu(z).\label{0314-4}
\end{eqnarray}
Letting $n\to\infty$, by (\ref{0314-1})-(\ref{0314-4}), we get
$$\|uC_\vp\|_{e,A_\om^p\to L_\mu^q}^q  \lesssim \sup_{|a|>r}\int_\D |F_{a,p,\gamma}(\vp(z))|^q |u(z)|^q d\mu(z)
+\varepsilon^q  \|u\|_{L^q_\mu}^q .$$
Since $\varepsilon$ is arbitrary, by letting $r\to 1$, we obtain
$$\|uC_\vp\|_{e,A_\om^p\to L_\mu^q}^q  \lesssim \limsup_{|a|\to 1}\int_\D |F_{a,p,\gamma}(\vp(z))|^q |u(z)|^q d\mu(z).$$

{\bf Lower estimate of $\|uC_\vp\|_{e,A_\om^p\to L_\mu^q}$. }

Assume that $K:A_\om^p\to L_\mu^q$ is compact. By Lemmas 2 and 5,   $\lim\limits_{|a|\to 1} \|KF_{a,p,\gamma}\|_{L_\mu^q}=0$.
Then
\begin{eqnarray*}
\|uC_\vp -K\|_{A_\om^p\to L_\mu^q}
&\gtrsim& \limsup_{|a|\to 1} \|(uC_\vp-K)F_{a,p,\gamma}\|_{L_\mu^q}   \\
&\geq& \limsup_{|a|\to 1} \|uC_\vp F_{a,p,\gamma}\|_{L_\mu^q}-\limsup_{|a|\to 1} \|K F_{a,p,\gamma}\|_{L_\mu^q} \\
&=& \limsup_{|a|\to 1} \|uC_\vp F_{a,p,\gamma}\|_{L_\mu^q}.
\end{eqnarray*}
Therefore, we get
$$\|uC_\vp\|_{e,A_\om^p\to L_\mu^q}^q  \gtrsim  \limsup_{|a|\to 1}\int_\D |F_{a,p,\gamma}(\vp(z))|^q |u(z)|^q d\mu(z),$$
as desired. The proof is complete.
\end{proof}

 {\noindent\bf Theorem 3.}  {\it Assume  $\om\in\hD, 0<q< p<\infty$, $u:\D\to\CC$ is a measurable function, $\vp$ is an analytic self-map of $\D$, and $\mu$ is a positive Borel measure on $\D$. Then the following statements are equivalent.
  \begin{enumerate}[(i)]
    \item $uC_\vp: A_\om^p\to L^q_\mu$ is bounded;
    \item $uC_\vp: A_\om^p\to L^q_\mu$ is compact;
    \item $\Psi_{u,\vp}^\gamma(a)\in L_\om^{\frac{p}{p-q}}$   for all $\gamma$ large enough;
    \item $\Psi_{u,\vp}^\gamma(a)\in L_\om^{\frac{p}{p-q}}$ for some $\gamma$  large enough.
  \end{enumerate}
  Moreover, if $\gamma$ is fixed,
  \begin{equation}\label{1104-1}
  \|uC_\vp\|^q _{A_\om^p\to L_\mu^q}\approx \|\Psi_\gamma\|_{L_\om^\frac{p}{p-q}}.
  \end{equation}
  Here,
  \begin{equation*}
  \Psi_{u,\vp}^\gamma(a):=\int_\D |F_{a,p,\gamma}(\vp(z))|^p  |u(z)|^q d\mu(z) .
  \end{equation*}
  }

\begin{proof} Let  $\upsilon(E)=\int_{\vp^{-1}(E)}|u(z)|^q d\mu(z)$ for all Borel set $E$. By (\ref{0313-1}), we have
$$
\|uC_\vp f\|_{ L_\mu^q}^q =\|f\|_{L_\upsilon^q}^q.
$$
So,   $uC_\vp:A_\om^p\to L_\mu^q$ is bounded (compact) if and only if  $Id:A_\om^p\to L_\upsilon^q$  is bounded (compact).
By (26) in \cite{Pja2015} and Theorem 3 in \cite{PjaRjSk2015mz}, we have (i)$\Leftrightarrow$(ii)$\Leftrightarrow$(iii) and (\ref{1104-1}).
Since $1-|a|\leq |1-\ol{a}z|$ holds for all $a,z\in\D$, we get (iii) $\Leftrightarrow$ (iv).
%
The proof is complete.
\end{proof}

  {\noindent\bf Remark 2.}   Suppose all of the $p,q,\mu,\om,u,\vp$ meet the conditions of the Theorem 3 and $\upsilon$ is defined as in the proof of Theorem 3.
  For all $a\in\D\backslash\{0
  \}$ and $r\in(0,1)$, let
 $$ \Gamma(a)=\left\{z\in\D: \left|\arg z-\arg a\right|<\frac{1}{2}\left(1-\frac{|z|}{|a|}\right)\right\}$$
 $$T(a)=\left\{z\in\D: a\in\Gamma(z)\right\},\,\,\Delta(a,r)=\left\{z\in\D:\left|\frac{a-z}{1-\ol{a}z}\right|<r\right\},$$
 and
 $$
 Q(z)=\int_\Gamma(z)\frac{d\upsilon(\xi)}{\om(T(\xi))},\,\,
 M_\om(\upsilon)(z)=\sup_{z\in S(a)}\frac{\upsilon(S(a))}{\om(S(a))},\,\,
 \Phi_r(z)=\int_{\Gamma(z)}\frac{\upsilon(\Delta (a,r))}{\om(T(z))}\frac{dA(z)}{(1-|z|)^2}.$$
 By \cite[Theorem 3]{PjaRjSk2015mz}, for any fixed $\gamma$ and $r$, we have
   \begin{align}\label{1104-2}
  \|uC_\vp\|^q _{A_\om^p\to L_\mu^q}
  \approx  \|Id\|^q _{A_\om^p\to L_\upsilon^q}
  \approx \|\Psi_\gamma\|_{L_\om^\frac{p}{p-q}}
  \approx \|M_\om(\upsilon)\|_{L_\om^\frac{p}{p-q}}
  \approx \|Q\|_{L_\om^\frac{p}{p-q}}
  \approx \|\Phi_r\|_{L_\om^\frac{p}{p-q}}.
  \end{align}
  Using (26) in \cite{Pja2015}, we know that all the $\om(S(a))$ and $\om(T(a))$ in (\ref{1104-1}) and (\ref{1104-2}) can be exchanged.
  \msk

 {\noindent\bf Theorem 4.}  {\it Assume  $\om\in\hD$ satisfying $\int_0^1 \left(\log\frac{e}{1-t}\right)^2\om(t)dt<\infty$, $0<p<\infty$, $u:\D\to\CC$ is a measurable function, $\vp$ is an analytic self-map of $\D$.  If $uC_\vp:A_\om^2\to A_\om^2$ is compact, then
 $uC_\vp\in S_p(A_\om^2)$ if and only if
$$\int_\D\left(\frac{\sigma(\Delta(z,r))}{\om_*(z)}\right)^\frac{p}{2}\frac{dA(z)}{(1-|z|^2)^2}<\infty
 $$
for some (equivalently for all) $0<r<1$. Moreover, we have
$$
  |uC_\vp|_p^p\approx \int_\D\left(\frac{\sigma(\Delta(z,r))}{\om_*(z)}\right)^\frac{p}{2}\frac{dA(z)}{(1-|z|^2)^2}.
$$
Here $\sigma(E)=\int_{\vp^{-1}(E)}|u(z)|^2 \om(z)dA(z)$ for all Borel set $E\subset\D$, and $|\cdot|_p$ is the norm of $p$-Schatten class of $A_\om^2$.}

\begin{proof} For all $f,g\in A_\om^2$, we have
\begin{align}
\langle (uC_\vp)^*uC_\vp f,g\rangle_{A_\om^2} &=\langle uC_\vp f,uC_\vp g\rangle_{A_\om^2}
=\int_\D f(\vp(z))\ol{g(\vp(z))}|u(z)|^2\om(z)dA(z)\nonumber\\
&=\int_\D f(\zeta)\ol{g(\zeta)}d\sigma(\zeta).  \label{1017-1}
\end{align}
Suppose $B_z^\om(\zeta)$ is the reproducing kernel of $A_\om^2$, that is,
$$f(z)=\langle f,B_z^\om \rangle_{A_\om^2}=\int_\D  f(\zeta)\ol{B_z^\om(\zeta)}\om(\zeta)dA(\zeta).$$
Consider the Toeplitz operator,
$$T_\sigma f(z)=\int_\D f(\eta)\ol{B_z^\om(\eta)}d\sigma(\eta).$$
Since $\om$ is radial, by \cite[section 4.1]{PjaRj2014book}, polynomials are dense in $A_\om^p$ for all $0<p<\infty$.
So if $f,g\in A_\om^2$, there are two  polynomial sequences $\{f_n\}_{n=1}^\infty$  and $\{g_n\}_{n=1}^\infty$ such that $$\lim\limits_{n\to\infty}\|f-f_n\|_{A_\om^2}=0, \mbox{  and } \lim\limits_{n\to\infty}\|g-g_n\|_{A_\om^2}=0.$$
 Since $uC_\vp:A_\om^2\to A_\om^2$ is compact, by Theorem 1 and Lemma 3, we have
 \begin{align*}
 \|uC_\vp\|^2_{A_\om^2\to A_\om^2}\approx
\sup_{a\in\D} \int_\D|F_{a,2,\gamma}(z)|^2d\sigma(z)
 \approx \sup_{a\in\D} \frac{\sigma(S(a))}{\om(S(a))}<\infty.
\end{align*}
Then,  $Id:A_\om^2\to A_\sigma^2$  and $Id:A_\om^1\to A_\sigma^1$ are bounded by  Lemma 3. So, we have $\lim\limits_{k\to\infty}\|g-g_k\|_{A_\sigma^2}=0$.

For any $h\in H^\infty\subset A_\om^1$, letting $M=\sup_{z\in\D} |h(z)|$, by \cite[Theorem C]{PjaRjSk2018jga} (see also \cite[Theorem 1]{PjaRj2016jmpa}),  there exist a $C=C(h,u,\vp,\om)$ such that
\begin{align*}
|T_\sigma h(z)| &\leq M\int_\D  |B_z^\om(\eta)|d\sigma(\eta) \leq M\|Id\|_{A_\om^1\to A_\sigma^1}\|B_z^\om\|_{A_\om^1}  \leq C \log \frac{e}{1-|z|}.
\end{align*}
Therefore,
  $$\|T_\sigma h\|_{A_\om^2}^2\lesssim \int_\D \left(\log \frac{e}{1-|z|}\right)^2\om(z)dA(z)<\infty.$$
 That is to say, $T_\sigma h\in A_\om^2$.

Since $g_n\in A_\om^2=(A_\om^2)^*$, for any $n\in\N$, by Lemma 11 of \cite{PjaRjSk2018jga},
\begin{align}\label{1020-1}
\langle T_\sigma f_n,g\rangle_{A_\om^2} =\lim_{k\to\infty} \langle T_\sigma f_n,g_k\rangle_{A_\om^2}
=\lim_{k\to\infty} \langle  f_n,g_k\rangle_{A_\sigma^2}=\int_\D  f_n(\eta)\ol{g(\eta)}d\sigma(\eta).
\end{align}
Since $g$ is  arbitrary, by (\ref{1017-1}) and (\ref{1020-1}), we have
\begin{align}\label{1020-2}
T_\sigma f_n=(uC_\vp)^*(uC_\vp)f_n.
\end{align}

For any fixed $z_0\in\D$, when $|z-z_0|<\frac{1+|z_0|}{2}$, by H\"older inequality and \cite[Lemma 6]{PjaRjSk2018jga}, we have
\begin{align*}
|T_\sigma f(z)-T_\sigma f_n(z)|
&\leq \int_\D |f(\eta)-f_n(\eta)||\ol{B_z^\om(\eta)}|d\sigma(\eta)  \\
&\lesssim \sup_{|z|<(|z_0|+1)/2}\frac{1}{\om(S(z))}\sqrt{\sigma(\D)} \|f-f_n\|_{A_\sigma^2}.
\end{align*}
So, we obtain $\lim\limits_{n\to\infty} T_\mu f_n(z)=T_\mu f(z)$ for all $z\in\D.$
Using (\ref{1020-2}), we have $$T_\sigma=(uC_\vp)^*uC_\vp.$$

By \cite[Theorem 1.26]{Zk2007book}, when $p>0$, $uC_\vp\in S_p(A_\om^2)$ if and only if $(uC_\vp)^*uC_\vp\in S_{\frac{p}{2}}(A_\om^2)$.
By \cite[Theorem 3]{PjaRjSk2018jga} (or \cite[Theorem 1]{PjaRj2016}), $T_\sigma\in S_\frac{p}{2}(A_\om^2)$ if and only if
$$\int_\D\left(\frac{\sigma(\Delta(z,r))}{\om_*(z)}\right)^\frac{p}{2}\frac{dA(z)}{(1-|z|^2)^2}<\infty,$$
for some (equivalently for all) $r\in (0,1)$. Moreover,
$$|uC_\vp|_p^p =|T_\sigma|_{\frac{p}{2}}^\frac{p}{2}
\approx \int_\D\left(\frac{\sigma(\Delta(z,r))}{\om_*(z)}\right)^\frac{p}{2}\frac{dA(z)}{(1-|z|^2)^2}.$$
The proof is complete.
\end{proof}

{\noindent\bf Theorem 5. } {\it Suppose $1<p<\infty$, $\om\in\R$, $\vp$ is a finite Blaschke product and $u\in A_\om^p$. If $uC_\vp$ is bounded on $A_\om^p$, then $u\in H^\infty(\D)$. }

\begin{proof}  By \cite[Corollary 7]{PjaRj2016jmpa}, we have $(A_\om^p)^*\simeq A_\om^{p^\p}$ where $\frac{1}{p}+\frac{1}{p^\p}=1$.
Let $B_z^\om(\zeta)$ be the reproducing kernel of $A_\om^2$.
By \cite[Lemma 6]{PjaRjSk2018jga}, $B_z^\om\in H^\infty\subset A_\om^{p^\p}$.
Since $uC_\vp$ is bounded on $A_\om^p$, so is $(uC_\vp)^*$ on $A_\om^{p^\p}$. For all $f\in A_\om^p$, we have,
\begin{align*}
\langle  (uC_\vp)^* B_z^\om, f\rangle_{A_\om^2}=\langle  B_z^\om, uC_\vp f\rangle_{A_\om^2}
=\overline{u(z)f(\vp(z))}=\overline{u(z)}\langle B_{\vp(z)}^\om,f\rangle_{A_\om^2}.
\end{align*}
Therefore,
\begin{equation*}
(uC_\vp)^* B_z^\om=\overline{u(z)}B_{\vp(z)}^\om.
\end{equation*}
So,
$$|u(z)|\|B_{\vp(z)}^\om\|_{A_\om^{p^\p}}=\|\overline{u(z)}B_{\vp(z)}^\om\|_{A_\om^{p^\p}}
=\|(uC_\vp)^* B_z^\om\|_{A_\om^{p^\p}}\leq \|(uC_\vp)^*\|_{A_\om^{p^\p}\to A_\om^{p^\p}}\|B_z^\om\|_{A_\om^{p^\p}}.
$$
Let $M=\|(uC_\vp)^*\|_{A_\om^{p^\p}\to A_\om^{p^\p}}$. By \cite[Theorem C]{PjaRjSk2018jga}, we have
\begin{equation}\label{0401-1}
|u(z)|^{p^\p}\leq M^{p^\p} \left(\frac{\|B_z^\om\|_{A_\om^{p^\p}}}{\|B_{\vp(z)}^\om\|_{A_\om^{p^\p}}}\right)^{p^\p}
\approx  M ^{p^\p}\left(\frac{\om(S(\vp(z)))}{\om(S(z))}\right)^{p^\p-1}.
\end{equation}

Suppose $\vp(z)=z^m\prod_{k=1}^n \frac{|a_k|}{a_k}\frac{a_k-z}{1-\overline{a_k}z}$.
Let
$$c=\max\{|a_k|:k=1,2,\cdots,n\}, \mbox{ and } d=\min\{|a_k|:k=1,2,\cdots,n\}.$$
As in the proof of \cite[Lemma 2.1]{CmdHag2003}, for $c< |z|<1$, we have
$$\frac{1-|\vp(z)|^2}{1-|z|^2}\leq m+2n\frac{1+d}{1-d}.$$

By Lemma 1, there are $1<a<b<\infty$ and $\delta\in(0,1)$, such that
$\frac{\om_*(r)}{(1-r)^b}$ is increasing on $ [\delta,1)$.
Let $$r_0=\inf\left\{r:r>\max\{c,\delta\} \mbox{ and }|\vp(z)|\geq \delta\mbox{ for all } |z|=r\right\}.$$
Then $0<r_0<1.$ Obviously, we have
\begin{align}\label{1017-2}
\sup_{|z|\leq r_0}\frac{\om(S(\vp(z)))}{\om(S(z))}<\infty.
\end{align}

When $|z|>r_0$,  by (\ref{1017-3}),   we have
$$\frac{\om(S(\vp(z)))}{\om(S(z))}
\approx \frac{\om_*(\vp(z))}{\om_*(z)}.$$
So, if $|\vp(z)|<|z|$,
\begin{equation}\label{0401-3}
\frac{\om_*(\vp(z))}{\om_*(z)}
=\frac{\frac{\om_*(\vp(z))}{(1-|\vp(z)|)^b}}{\frac{\om_*(z)}{(1-|z|)^b}}  \frac{(1-|\vp(z)|)^b}{(1-|z|)^b}
\lesssim\left(m+2n\frac{1+d}{1-d}\right)^b.
\end{equation}
If $|z|\leq |\vp(z)|< 1$, by Lemma 1, we get
\begin{equation}\label{0401-2}
\frac{\om_*(\vp(z))}{\om_*(z)}
\approx
\frac{(1-|\vp(z)|)\int_{|\vp(z)|}^1 \om(t)dt}{(1-|z|)\int_{|z|}^1 \om(t)dt}
\lesssim
m+2n\frac{1+d}{1-d}.
\end{equation}
Therefore, by (\ref{0401-1})-(\ref{0401-2}), we obtain $u\in H^\infty.$ The proof is complete.
\end{proof}

To state and prove the next result, we introduce a set. Let $1<p<\infty$, $\om\in\R$ and $\vp$ be an analytic self-map of $\D$. We define
$$X:=\{u\in H(\D):uC_\vp(A_\om^p)\subset A_\om^p\}.$$

{\noindent\bf Theorem 6. } {\it Let $1<p<\infty$ and  $\vp$ be an analytic self-map of $\D$. Suppose $\om\in\R$ such that
\begin{enumerate}[(i)]
  \item $\hat\om(\vp_t(r))\hat\om(r)\lesssim \hat\om(t),\mbox{ for all } 0\leq r\leq t<1$, here $\vp_t(r)=\frac{t-r}{1-tr}$;
  \item $2A+AB-B> 0$, where $A=\liminf\limits_{t\to 1}\frac{\int_t^1 \om(s)ds}{(1-t)\om(t)} $ and $B=\limsup\limits_{t\to 1}\frac{\int_t^1 \om(s)ds}{(1-t)\om(t)} $.
\end{enumerate}
If $X =H^\infty$, then $\vp$ is a finite Blaschke product. }

\begin{proof}  By \cite[Proposition 18]{PjaRj2016}, for any $\om\in\hD$ and $\vp$ is an analytic self map of $\D$,  $C_\vp:A_\om^p\to A_\om^p$ is bounded.
So, for  any $u\in X$, we can define $\|u\|_X=\|uC_\vp\|_{A_\om^p\to A_\om^p}$. Let $\{u_n\}$ be a Cauchy sequence in $X$. Then $\{u_nC_\vp\}$ is a Cauchy sequence in $B(A_\om^p)$.
Here, $B(A_\om^p)$ means the set of bounded operators on   $A_\om^p$.
 So, there exists a $T\in B(A_\om^p)$, such that $\lim\limits_{n\to\infty} u_n C_\vp=T$.
Since $f(z)=1\in A_\om^p$,
$$u:=T(1)\in A_\om^p, \mbox{ and } \lim\limits_{n\to\infty}\|u_n-u\|_{A_\om^p}=0.$$
Therefore, for all $f\in A_\om^p$,
$$\lim\limits_{n\to\infty} u_n(z)f(\vp(z))=u(z)f(\vp(z)).$$
Since $\lim\limits_{n\to\infty}\|u_nC_\vp f-Tf\|_{A_\om^p}=0$,
$$\lim\limits_{n\to\infty} u_n f(\vp(z))=(Tf)(z).$$
So, we have $Tf=uC_\vp f$ for all $f\in A_\om^p$. Therefore $u\in X$. That is to say, $X$ is complete under the norm $\|\cdot\|_X$.

Since $X=H^\infty$, $C_\vp\in B(A_\om^p)$. Then for all $u\in X$, $\|u\|_X\leq \|u\|_{H^\infty} \|C_\vp\|_{A_\om^p\to A_\om^p}$.
By Inverse Mapping Theorem, $\|u\|_X\approx \|u\|_{H^\infty}$.

By $\om\in \R$, we have $AB>0$. Therefore $\frac{2}{B}+1> \frac{1}{A}$.
So, there exists $\varepsilon\in (0,A)$ such that $\frac{2}{B+\varepsilon}+1-\frac{1}{A-\varepsilon}>0.$

Let $a=\frac{1}{B+\varepsilon}-1$ and $b=\frac{1}{A-\varepsilon}-1$. Then $2a+2-b>0.$
By the proof of  \cite[Lemma 1.1]{PjaRj2014book}, there is a $\delta\in(0,1)$ such that
$$
\frac{\om (r)}{(1-r)^a} \;\; \mbox{is essential decreasing on}\;\; [\delta,1)\;
\mbox{and}\;\; \lim_{r\to 1}\frac{\om (r)}{(1-r)^a}=0;     $$
and
$$
\frac{\om (r)}{(1-r)^b} \;\; \mbox{is essential increasing on}\;\; [\delta,1)\;\;
\mbox{and}\;\; \lim_{r\to 1}\frac{\om (r)}{(1-r)^b}=\infty.
$$
Let
$$
\mu(z)=\left\{
\begin{array}{cc}
  \om(z), & \delta \leq |z|<1; \\
  \frac{\om(\delta)(1-|z|)^a}{(1-\delta)^a}, & |z|<\delta.
\end{array}
\right.
$$
Then it is easy to check that the following statements hold.
\begin{enumerate}[(i)]
  \item $\frac{\mu (r)}{(1-r)^a} \;\; \mbox{is essential decreasing on}\;\; [0,1)\;
\mbox{and}\;\; \lim_{r\to 1}\frac{\mu (r)}{(1-r)^a}=0;     $
  \item $
\frac{\mu (r)}{(1-r)^b} \;\; \mbox{is essential increasing on}\;\; [0,1)\;\;
\mbox{and}\;\; \lim_{r\to 1}\frac{\mu (r)}{(1-r)^b}=\infty;
$
  \item $\mu\in\R$,  $A=\liminf\limits_{t\to 1}\frac{\int_t^1 \mu(s)ds}{(1-t)\mu(t)} $ and $B=\limsup\limits_{t\to 1}\frac{\int_t^1 \mu(s)ds}{(1-t)\mu(t)}; $
  \item $\|f\|_{A_\mu^p}\approx \|f\|_{A_\om^p}$ and $\hat\mu(\vp_t(r))\hat\mu(r)\lesssim \hat\mu(t),\mbox{ for all } 0\leq r\leq t<1$.
\end{enumerate}
Therefore, without loss of generality, let $\delta=0$.  So,  we have
\begin{align*}
\frac{\om(\vp_w(z))}{\om(z)}
=\frac{\frac{\om(\vp_w(z))}{1-|\vp_w(z)|^a}(1-|\vp_w(z)|)^a}{\frac{\om(z)}{(1-|z|)^a}(1-|z|)^a}
&\lesssim\left(\frac{1-|\vp_w(z)|^2}{1-|z|^2}\right)^a,
\mbox{ when } |\vp(z)|> |z|,
\end{align*}
and
\begin{align*}
\frac{\om(\vp_w(z))}{\om(z)}
=\frac{\frac{\om(\vp_w(z))}{1-|\vp_w(z)|^b}(1-|\vp_w(z)|)^b}{\frac{\om(z)}{(1-|z|)^b}(1-|z|)^b}
&\lesssim\left(\frac{1-|\vp_w(z)|^2}{1-|z|^2}\right)^b,
\mbox{ when } |\vp(z)|\leq |z|.
\end{align*}
Therefore,
\begin{align}\label{0521-3}
\frac{\om(\vp_w(z))}{\om(z)}
&\lesssim\left(\frac{1-|\vp_w(z)|^2}{1-|z|^2}\right)^a+\left(\frac{1-|\vp_w(z)|^2}{1-|z|^2}\right)^b
= \left(\frac{1-|w|^2}{|1-\ol{w}z|^2}\right)^a+\left(\frac{1-|w|^2}{|1-\ol{w}z|^2}\right)^b.
\end{align}
Let $\alpha=  \frac{2(a+2)}{p}$ and  $u_w(z)=\left(\frac{1}{1-\ol{w}z}\right)^\alpha$. Then $\|u_w\|_\infty^p=\frac{1}{(1-|w|)^{2a+4}}$. For all $f\in A_\om^p$, by (\ref{0521-3}), we have
\begin{align*}
\|u_wC_\vp f\|_{A_\om^p}^p
&=\int_\D \frac{1}{|1-\ol{w}z|^{p\alpha}}|f\circ \vp(z)|^p\om(z)dA(z)\\
&=\int_\D \frac{1}{|1-\ol{w}\vp_w(z)|^{p\alpha}}|f\circ\vp\circ \vp_w(z)|^p  |\vp_w^\p(z)|^2\om(\vp_w(z)) dA(z)   \\
&\lesssim\frac{1}{(1-|w|^2)^{a+2}}  \int_\D  |f\circ\vp\circ \vp_w(z)|^p
      \left(1+    \left(\frac{1-|w|^2}{|1-\ol{w}z|^2}\right)^{b-a}\right)\om(z)dA(z)\\
&\lesssim \frac{1}{(1-|w|^2)^{b+2}} \int_\D  |f\circ\vp\circ \vp_w(z)|^p  \om(z)dA(z)\\
&\lesssim \frac{\|f\|_{A_\om^p}^p}{(1-|w|^2)^{b+2}(1-|\vp(w)|)\hat{\om}(\vp(w))}.
\end{align*}
The last inequality can be got by  \cite[Proposition 18]{PjaRj2016}. By $\|u_w\|_\infty\approx \|u_wC_\vp\|_{A_\om^p\to A_\om
^p}$, we have
$$\frac{1}{(1-|w|^2)^{2a+2-b}}\lesssim \frac{1}{(1-|\vp(w)|)\hat{\om}(\vp(w))}.  $$
By $2a+2-b>0$,  we have $|\vp(w)|\to 1$ as $|w|\to 1$.
By Lemma 6, $\vp$ is a finite Blaschke product. The proof is complete.
\end{proof}

By Theorems 5 and 6, for some regular weight $\om$, $X= H^\infty$  if and only if $\vp$ is a finite Blaschke product. Here, we give two examples.\msk

{\noindent\bf Corollary 7. }{\it Let $1<p<\infty$, and  $\vp$ be an analytic self map of $\D$. Suppose $\om$ is  either (a) or (b).
 \begin{enumerate}[(a)]
   \item $\om(r)=(1-r)^\alpha\left(\log \frac{e}{1-r}\right)^\beta,$    $\alpha>-1$ and $\beta\leq 0$;
   \item $\om(r)=\exp\left(-\beta\left(\log\frac{e}{1-r}\right)^\alpha\right)$,   $0<\alpha\leq 1$ and $\beta>0$.
 \end{enumerate}
Then, $X=H^\infty$ if and only if $\vp$ is a finite Blaschke product.}

  \begin{proof}By \cite[(4.4) and (4.5)]{AaSa1997iuaj}, the weights (a) and (b)  are regular.

Suppose  $\vp$ is a finite Blaschke product. By \cite[Theorem 18]{PjaRj2016}, $C_\vp:A_\om^p\to A_\om^p$ is bounded. So,  $H^\infty\subset X$.
By Theorem 5, $X\subset H^\infty.$ Therefore, $X=H^\infty.$

Suppose $X=H^\infty.$ By Bernouilli-l'H$\hat{\mathrm{o}}$pital theorem,  both  (a) and (b) meet the condition (ii) of Theorem 6.
  So, if we can prove that (a) and (b) meet the condition (i) of Theorem 6, then $\vp$ is a finite Blaschke product.

  {\bf Condition (a).} When $0\leq r\leq t<1$,  let $\theta=r/t$ and
    $$f(\theta,t)=\log\frac{e(1-\theta t^2)}{(1-t)(1+\theta t)}\log\frac{e}{1-\theta t}, \,\,(0\leq \theta\leq 1, 0< t<1).$$
Then
$$f_\theta^\p(\theta,t)=-\left(\frac{t^2}{1-\theta t^2}+\frac{t}{1+\theta t}\right)\log\frac{e}{1-\theta t}+
\frac{t}{1-\theta t}\log\frac{e(1-\theta t^2)}{(1-t)(1+\theta t)}.$$
Suppose $t>0$, and let
$$g(\theta,t)=\frac{1-\theta t}{t}f_\theta^\p(\theta,t)=\log\frac{e(1-\theta t^2)}{(1-t)(1+\theta t)}
-\left(\frac{t(1-\theta t)}{1-\theta t^2}+\frac{1-\theta t}{1+\theta t}\right)\log\frac{e}{1-\theta t}.
 $$
Then
\begin{align*}
g_\theta^\p(\theta,t) &=-\frac{2t^2}{1-\theta t^2}-\frac{2t}{1+\theta t}-
\left(\frac{t^3-t^2}{(1-\theta t^2)^2}-\frac{2t}{(1+\theta t)^2}\right)\log\frac{e}{1-\theta t}\\
&=\frac{1}{(1+\theta t)^2}\left(   \frac{(1+\theta t)^2}{(1-\theta t^2)^2}h(\theta,t)+2tk(\theta,t)      \right),
\end{align*}
where
$$h(\theta,t)=(t^2-t^3)\log\frac{e}{1-\theta t}-2t^2(1-\theta t^2)$$
and
$$k(\theta,t)=\log\frac{e}{1-\theta t}-(1+\theta t).$$
Since
$$h_\theta^\p(\theta,t)=\frac{t^3+t^4-2\theta t^5}{1-\theta t}>0, \,\, k_\theta^\p(\theta,t)=t\left(-1+\frac{1}{1-\theta t}\right)>0,$$
we have
 $$G(\theta,t):= (1+\theta t)^2  g_\theta^\p(\theta,t)$$
is increasing on $[0,1]$ about $\theta$.
Since $\lim\limits_{t\to 1}G(1,t)=+\infty$, there exists a $\tau\in (0,1)$ such that
$G(1,t)>0$,  for all  $t\in (\tau,1).$
If $t\in(\tau,1)$, by $G(0,t)<0$,
there is a $\nu(t)\in (0,1)$, such that
$$G(\theta,t)<0,  \mbox{ when }\theta\in [0,\nu(t)) ,$$
and
$$G(\theta,t)>0,   \mbox{ when }\theta\in (\nu(t),1].$$
Since $g^\p_\theta(\theta,t)=\frac{G(\theta,t)}{(1+\theta t)^2}$, when $t\in (\tau,1)$,  $g(\theta,t)$ is decreasing on $ [0,\nu(t))$ and increasing on $(\nu(t),1]$.
Since
$$g(0,t)=\log\frac{e}{1-t}-(t+1)>0,\,\mbox{and}\,g(1,t)=\frac{1}{t+1}\left(t+1-\log\frac{e}{1-t}\right)<0,$$
 So there is a $\mu(t)\in (0,1)$ for every $t\in(\tau,1)$,
 such that,   $f(\theta,t)$ is increasing on $[0,\mu(t))$ and decreasing on $(\mu(t),1]$.
 Since  $f(0,t)=f(1,t)=\log\frac{e}{1-t}$,
 $$\frac{f(\theta,t)}{\log\frac{e}{1-t}}\geq 1, \mbox{ when } t\in(\tau,1) \mbox{ and } \theta\in[0,1].$$
 It is obvious that
$$\inf_{t\in [0,\tau] ,{\theta\in [0,1]}}\frac{f(\theta,t)}{\log\frac{e}{1-t}}>0.$$
Therefore, we have
$$C_1:=\inf_{0\leq r\leq t<1}\frac{\log\frac{e(1-rt)}{(1-t)(1+r)}\log\frac{e}{(1-r)}}{\log\frac{e}{1-t}}>0.$$
So, when $\alpha>-1,\beta\leq 0$ and $\om(r)=(1-r)^\alpha\left(\log \frac{e}{1-r}\right)^\beta$, we have
\begin{align*}
\frac{{\om(\vp_t(r))\om(r)}}{ {\om(t)}}
\approx \left(\frac{\log\frac{e(1-rt)}{(1-t)(1+r)}\log\frac{e}{(1-r)}}{\log\frac{e}{1-t}}\right)^\beta \leq {C_1^\beta}.
\end{align*}
Since $\om\in\R$,   we get
\begin{align}\label{1012-1}
\frac {{\hat\om(\vp_t(r))\hat\om(r)}  }  { {\hat\om(t)}}   \approx \frac{  {\om(\vp_t(r))\om(r)}}{ {\om(t)}}.
\end{align}
Therefore,
$$ {{\hat\om(\vp_t(r))\hat\om(r)}  } \lesssim { {\hat\om(t)}}, \mbox{ when } 0\leq r\leq t<1.$$

{\bf Condition (b).} Suppose $0\leq r\leq t<1$. Since
$\frac{e(1-rt)}{1-r^2}\geq \frac{e(1-r)}{1-r^2}>1$, when $0<\alpha\leq 1$, we have
\begin{align*}
&\left(\log\frac{e(1-rt)}{(1-t)(1+r)}\right)^\alpha + \left(\log\frac{e}{(1-r)}\right)^\alpha
\geq    \left(\log \frac{e^2(1-rt)}{(1-t)(1-r^2)}\right)^\alpha
\geq \left(\log\frac{e}{1-t}\right)^\alpha.
\end{align*}
So, when $0<\alpha\leq 1$, $\beta>0$ and $\om(r)=\exp\left(-\beta\left(\log\frac{e}{1-r}\right)^\alpha\right)$, we have
 $${\om(\vp_t(r))\om(r)}\lesssim {\om(t)}.$$
By (\ref{1012-1}), we get
$$ {{\hat\om(\vp_t(r))\hat\om(r)}  } \lesssim { {\hat\om(t)}}, \mbox{ when } 0\leq r\leq t<1.$$
 The proof is complete.
\end{proof}

\end{document}